\documentclass[12pt]{article}
\usepackage{latexsym,amsmath,amssymb}

\newfont{\bb}{msbm10 at 12pt}

\def\r{\hbox{\bb R}}

\def\h{\hbox{\bb H}}

\def\c{\hbox{\bb C}}
\def\s{\hbox{\bb S}}
\def\mr{\mathbb{M}^2\times\mathbb{R}}
\def\hr{\mathbb{H}^2\times\mathbb{R}}
\def\mkr{\mathbb{M}^2(\kappa)\times\mathbb{R}}
\def\sr{\mathbb{S}^2\times\mathbb{R}}
\def\m{\mathbb{M}^2}
\def\mk{\mathbb{M}^2(\kappa)}

\newcommand{\norm}[1]{\left\Vert #1 \right\Vert}

\newcommand{\set}[1]{\left\{#1\right\}}
\newcommand{\meta}[2]{\langle #1,#2 \rangle }

\newcommand{\camb}{\overline{\nabla}}
\newcommand{\ext}{\wedge}
\newcommand{\To}{\longrightarrow }

\newcommand{\parz}{\partial_z}

\newcommand{\zb}{\bar{z}}
\newcommand{\hm}{\hbox{\bb E}(\kappa , \tau)}
\newcommand{\Ab}{\overline{A}}
\newcommand{\pb}{\overline{p}}

\usepackage[latin1]{inputenc}
\usepackage{amsmath}
\usepackage{amsthm}
\usepackage{times}
\usepackage{amscd}
\usepackage{epsf}

\numberwithin{equation}{section}

\begin{document}

\theoremstyle{plain}\newtheorem{lemma}{Lemma}[section]
\theoremstyle{plain}\newtheorem{proposition}{Proposition}[section]
\theoremstyle{plain}\newtheorem{theorem}{Theorem}[section]
\theoremstyle{plain}\newtheorem{definition}{Definition}[section]
\theoremstyle{plain}\newtheorem{remark}{Remark}[section]
\theoremstyle{plain}\newtheorem{corollary}{Corollary}[section]

\begin{center}
\rule{15cm}{1.5pt} \vspace{.4cm}

{\Large \bf Complete Constant Mean Curvature surfaces \\[4mm] in homogeneous spaces}
\vspace{0.5cm}

{\large Jos$\acute{\text{e}}$ M. Espinar$\,^\dag$\footnote{The author is partially
supported by Spanish MEC-FEDER Grant MTM2010-19821, and Regional J. Andalucia Grants
P06-FQM-01642 and FQM325}, Harold Rosenberg$\,^\ddag$}\\
\vspace{0.4cm}\rule{15cm}{1.5pt}
\end{center}

\vspace{.4cm}

\noindent $\mbox{}^\dag$ Institut de Mathématiques, Universit$\acute{\text{e}}$
Paris VII, 175 Rue du Chevaleret, 75013 Paris, France; e-mail:
jespinar@ugr.es\vspace{0.2cm}

\noindent $\mbox{}^\ddag$ Instituto de Matematica Pura y Aplicada, 110 Estrada Dona
Castorina, Rio de Janeiro 22460-320, Brazil; e-mail: rosen@impa.br

\vspace{.3cm}

\begin{abstract}
In this paper we classify complete surfaces of constant mean curvature whose
Gaussian curvature does not change sign in a simply connected homogeneous manifold
with a 4-dimensional isometry group.
\end{abstract}

\section{Introduction}

%One of the major efforts of mathematicians in Differential Geometry is the
%classification of complete constant mean curvature surfaces, in short $H-$surfaces.
%Despite what happens for complete constant Gaussian curvature surfaces in $\r ^3$,
%such a classification was finished in the earliest 1950's, the classification of
%complete $H-$surfaces in $\r ^3$ with no extra hypothesis is impossible. In this
%line for $H\neq 0$, under topological hypothesis, we have Hopf Theorem \cite{Ho},
%any topological sphere with constant mean curvature in $\r ^3$ is a round sphere, or
%Nitsche Theorem \cite{N}, any topological disk of constant mean curvature in $\r^3$,
%whose boundary is a line of curvature, is a part of either a round sphere or a
%plane. Under {\it injectivity} conditions, Alexandrov Theorem \cite{A} states that
%the only compact embedded $H-$surfaces in $\r ^3$ are the round spheres, and
%Korevaar-Kusner-Solomon Theorem \cite{KKS} classify properly embedded $H-$surfaces
%with finite topology and less than 3 ends. For minimal surfaces, $H=0$, this task
%seems even harder. {\huge Meeks-Rosenberg --- Meeks-P\'{e}rez-Ros.... Write
%something here!}

In 1966, T. Klotz and R. Ossermann showed the following:

\vspace{.2cm}

{\bf Theorem \cite{KO}:} {\it A complete $H-$surface in $\r ^3 $ whose Gaussian
curvature $K$ does not change sign is either a sphere, a minimal surface, or a right
circular cylinder.}

\vspace{.2cm}

The above result was extended to $\s^3 $ by D. Hoffman \cite{H}, and to $\h ^3$ by
R. Tribuzy \cite{T} with an extra hypothesis if $K$ is non-positive. The additional
hypothesis says that, when $K\leq 0$, one has $H^2 - K -1 >0$.

In recent years, the study of $H-$surfaces in product spaces and, more generally, in
a homogeneous three-manifold with a 4-dimensional isometry group is quite active
(see \cite{AR,AR2}, \cite{CoR}, \cite{ER}, \cite{FM,FM2}, \cite{DH} and references
therein).

The aim of this paper is to extend the above Theorem to homogeneous spaces with a
4-dimensional isometry group. These homogeneous space are denoted by $\hm$, where
$\kappa$ and $\tau$ are constant and $\kappa -4\tau ^2 \neq 0$. They can be
classified as $\mkr$ if $\tau = 0$, with $\mk = \s ^2(\kappa)$ if $\kappa >0$ ($\s
^2(\kappa)$ the sphere of curvature $\kappa$), and  $\mk  =  \h ^2(\kappa)$ if
$\kappa < 0$ ($\h^2(\kappa)$ the hyperbolic plane of curvature $\kappa$).  If $\tau$
is not equal to zero, $\hm $ is a Berger sphere if $\kappa > 0$, a Heisenberg space
if $\kappa = 0$ (of bundle curvature $\tau$), and the universal cover of ${\rm
PSL}(2,\r)$ if $\kappa < 0$. Henceforth we will suppose $\kappa$ is plus or minus
one or zero.

%They can be classified as: the product spaces $\hr$ if $\kappa = -1$ and $\tau =0$,
%or $\sr$ if $\kappa = 1$ and $\tau =0$, the Heisenberg space ${\rm Nil}_3$ if
%$\kappa = 0$ and $\tau = 1/2$, the Berger spheres $\s ^3 _{Berger}$ if $\kappa =1 $
%and $\tau \neq 0$, and the universal covering of ${\rm PSL}(2,\r)$ if $\kappa = -1$
%and $\tau \neq 0$.

The paper is organized as follows. In Section \ref{Sgeometry}, we establish the
definitions and necessary equations for an $H-$surface. We also state here two
classification results for $H-$surfaces. We prove them in Section \ref{Snuconstant}
and Section \ref{Sqconstant} for the sake of completeness.

Section \ref{SKnonnegative} is devoted to the classification of $H-$surfaces with
non-negative Gaussian curvature,

\vspace{.2cm}

{\bf Theorem \ref{Kgeq0}.} {\it Let $\Sigma \subset \hm$ be a complete $H-$surface
with $K\geq 0$. Then, $\Sigma$ is either a rotational sphere (in particular,
$4H^2+\kappa >0$), or a complete vertical cylinder over a complete curve of geodesic
curvature $2H$ on $\m (\kappa)$.}

\vspace{.2cm}

In Section \ref{SKnonpositive} we continue with the classification of $H-$surfaces
with non-positive Gaussian curvature.

\vspace{.2cm}

{\bf Theorem \ref{Kleq0}.} {\it Let $\Sigma \subset \hm$ be a complete $H-$surface
with $K\leq 0$ and $H^2+\tau^2 - \left|\kappa -4\tau ^2\right|
>0$. Then, $\Sigma$ is a complete vertical cylinder over a complete curve of
geodesic curvature $2H$ on $\m (\kappa)$.}

\vspace{.2cm}

The above theorem is not true without the inequality; for example, any complete
minimal surface in $\hr$ that is not a vertical cylinder.

In the Appendix, we give a result, which we think is of independent interest,
concerning differential operators on a Riemannian surface $\Sigma$ of the form
$\Delta + g$, acting on $C^2(\Sigma)-$functions, where $\Delta $ is the Laplacian
with respect to the Riemannian metric on $\Sigma$ and $g \in C^0 (\Sigma)$.

\section{The geometry of surfaces in homogeneous spaces}\label{Sgeometry}

Henceforth $\hm $ denotes a complete simply connected homogeneous three-manifold
with $4-$dimensional isometry group. Such a three-manifold can be classified in
terms of a pair of real numbers $(\kappa , \tau)$ satisfying $\kappa - 4 \tau ^2
\neq 0$. In fact, these manifolds are Riemannian submersions over a complete
simply-connected surface $\m (\kappa)$ of constant curvature $\kappa$, $\pi : \hm
\To \m (\kappa )$, and translations along the fibers are isometries, therefore they
generate a Killing field $\xi$, called the \emph{vertical field}. Moreover, $\tau$
is the real number such that $\camb _X \xi = \tau X \ext \xi$ for all vector fields
$X$ on the manifold. Here, $\camb$ is the Levi-Civita connection of the manifold and
$\ext$ is the cross product.

%
%surface with $\partial \m = \emptyset$. Let $g$ be the metric of $\m $ and $\nabla $
%the Levi-Civita connection on $\mr $ with the product metric $\meta{}{}=g+dt^2$.

Let $\Sigma $ be a complete $H-$surface immersed in $\hm$. By passing to a
$2-$sheeted covering space of $\Sigma$, we can assume $\Sigma$ is orientable. Let
$N$ be a unit normal to $\Sigma$. In terms of a conformal parameter $z$ of $\Sigma$,
the first, $\meta{ \cdot }{ \cdot }$, and second, $II$, fundamental forms are given
by
\begin{equation}\label{I}
\begin{array}{l}
\meta{ \cdot }{ \cdot }\ =\ \lambda\,|dz|^2\\
II\ =\ p\,dz^2+\lambda\,H\,|dz|^2+\overline{p}\,d\bar{z}^2,
\end{array}
\end{equation}
where $p\,dz^2=\meta{ -\nabla _{\parz} N}{\parz}\,dz^2$ is the Hopf differential of
$\Sigma$.

Set $\nu = \meta{N}{\xi}$ and $T = \xi - \nu N$, i.e., $\nu $ is the normal
component of the vertical field $\xi$, called the \emph{angle function}, and $T$ is
the tangent component of the vertical field.

First we state the following necessary equations on $\Sigma$ which were obtained in
\cite{FM}.

\begin{lemma}\label{l1}
Given an immersed surface $\Sigma \subset \hm $, the following equations are
satisfied:
\begin{eqnarray}
K &=& K_e + \tau ^2 + (\kappa -4\tau ^2 )\, \nu ^2 \label{gauss}\\
p_{\zb}&=&\frac{\lambda}{2}\,(H_z+(\kappa -4\tau ^2)\,\nu\,A)\label{c1}\\
A_{\zb}&=&\frac{\lambda}{2}\,(H + i \tau )\,\nu\label{c2}\\
\nu_z&=&-(H-i\tau)\,A-\frac{2}{\lambda}\,p\,\Ab\label{c3}\\
|A|^2&=&\frac{1}{4}\, \lambda\,(1-\nu^2)\label{c4}\\
A_z &=&\frac{\lambda_z}{\lambda}\,A +p\,\nu\label{c5}
\end{eqnarray}
where $A= \meta{\xi}{\parz}$, $K_e$ the extrinsic curvature and $K$ the Gauss
curvature of $\Sigma$.
\end{lemma}

For an immersed $H-$surface $\Sigma \subset \hm$ there is a globally defined
quadratic differential, called the \emph{Abresch-Rosenberg} differential, which in
these coordinates is given by (see \cite{AR2}):

$$ Q \, dz^2 = (2(H +i\tau )\, p - (\kappa -4\tau ^2) A ^2)\,dz^2
,$$following the notation above.

It is not hard to verify this quadratic differential is holomorphic on an
$H-$surface using \eqref{c1} and \eqref{c2},

\begin{theorem}[\cite{AR},\cite{AR2}]\label{ARth}
$Q\, dz^2$ is a holomorphic quadratic differential on any $H-$surface in $\hm$.
\end{theorem}

Associated to the Abresch-Rosenberg differential we define the smooth function $q :
\Sigma \To [0 , + \infty )$ given by
$$ q = \frac{4 |Q|^2}{\lambda ^2} .$$
By means of Theorem \ref{ARth}, $q$ either has isolated zeroes or vanishes
identically. Note that $q$ does not depend on the conformal parameter $z$, hence $q$
is globally defined on $\Sigma$.

We continue this Section establishing some formulae relating the angle function, $q$
and the Gaussian curvature.

\begin{lemma}\label{l2}
Let $\Sigma$ be an $H-$surface immersed in $ \hm $. Then the following equations are
satisfied:
\begin{equation}\label{modnuz}
\begin{split}
\norm{\nabla \nu }^2 &= \dfrac{4H^2 +\kappa - (\kappa -4\tau ^2)\nu ^2}{4(\kappa
-4\tau ^2)}\left(4(H^2-K_e)+(\kappa -4\tau ^2)(1-\nu^2) \right)- \dfrac{q}{\kappa
-4\tau ^2}
%\\[3mm]
% & \, \, \left( 4H^2 +\kappa -(\kappa -4\tau ^2)\nu ^2\right)\frac{H^2-K_e}{\kappa -4\tau ^2}- \dfrac{q}{\kappa -4\tau ^2} .
\end{split}
\end{equation}

\begin{equation}\label{deltanu}
\Delta \nu = -\left( 4H^2 + 2\tau ^2 +(\kappa -4\tau ^2)(1-\nu ^2)-2K_e\right)\nu .
\end{equation}

Moreover, away from the isolated zeroes of $q$, we have
\begin{equation}\label{lapGauss}
\Delta \ln q = 4 K .
\end{equation}
\end{lemma}
\begin{proof}
From \eqref{c3}
$$
|\nu_z|^2\ =\ \frac{4\,|p|^2\,|A|^2}{\lambda^2}+ (H^2 + \tau
^2)\,|A|^2+\frac{2\,(H+i \tau)}{\lambda}p \Ab ^2 +\frac{2\,(H-i \tau)}{\lambda}\pb A
^2,$$ and taking into account that
$$ |Q |^2\ =\
4\,(H^2+\tau ^2)\,|p|^2+(\kappa - 4 \tau ^2)^2|A|^4-(\kappa - 4 \tau
^2)\left(2\,(H+i \tau)p \Ab ^2 +2\,(H-i \tau)\pb A ^2\right),
$$
we obtain, using also \eqref{c4}, that
\begin{equation*}
\begin{split}
|\nu_z|^2&= (H^2 + \tau ^2)|A|^2 + (H^2 -K_e)|A|^2 + (\kappa -4\tau
^2)\frac{|A|^4}{\lambda} \\
 & + 4 \left( \dfrac{H^2 +\tau ^2}{\kappa - 4\tau ^2}\right)\frac{|p|^2}{\lambda} -
 \dfrac{|Q|^2}{(\kappa -4\tau ^2)\lambda}
\end{split}
\end{equation*}where we have used that $4|p|^2= \lambda ^2 (H^2 -K_e)$ and $\kappa -4\tau^2 \neq 0$. Thus
\begin{equation*}
\begin{split}
\norm{\nabla \nu}^2 = \frac{4}{\lambda } |\nu _z|^2 &=(2H^2 -K_e + \tau ^2)(1-\nu^2)
+\dfrac{\kappa -4\tau^2}{4}(1-\nu^2)^2 \\
 & + 4 \left( \dfrac{H^2 +\tau ^2}{\kappa - 4\tau ^2}\right)(H^2 -K_e) -
 \dfrac{q}{\kappa -4\tau ^2} ,
\end{split}
\end{equation*}finally, re-ordering in terms of $H^2-K_e$ we have the expression.

On the other hand, by differentiating \eqref{c3} with respect to $\zb$ and using
\eqref{c5}, \eqref{c2} and \eqref{c1}, one gets
\begin{equation*}
\nu_{z\zb}\ = -(\kappa -4 \tau ^2
)\,\nu\,|A|^2-\,\frac{2}{\lambda}\,|p|^2\,\nu-\,\frac{H^2+\tau ^2}{2}\,\lambda\,\nu.
\end{equation*}

Then, from \eqref{c4},
\begin{equation*}
\nu_{z\zb}=-\frac{\lambda\,\nu}{4}\left((\kappa - 4\tau ^2)(1-\nu^2)
+\,\frac{8\,|p|^2}{\lambda^2}+2\,(H^2+\tau ^2)\right)
\end{equation*}thus
\begin{equation*}
\begin{split}
\Delta \nu &= \frac{4}{\lambda}\nu _{z\zb} =-\left( (\kappa - 4\tau ^2)(1-\nu^2)
+2(H^2 -K_e)+2\,(H^2+\tau ^2)\right)\nu .
\end{split}
\end{equation*}

Finally,
$$\Delta \ln q = \Delta \ln \dfrac{4 |Q|^2}{\lambda ^2 } = - 2 \Delta \ln \lambda = 4 K ,
$$where we have used that $Q\, dz ^2 $ is holomorphic and the expression of the
Gaussian curvature in terms of a conformal parameter.
\end{proof}

\begin{remark}
Note that \eqref{deltanu} is nothing but the Jacobi equation for the Jacobi field
$\nu$.
\end{remark}

Next, we recall a definition in these homogeneous spaces.

%(see \cite{LH}). Let $\gamma _s$ denote the integral curves of $\xi$. Consider a
%smooth embedding $\varphi: \m (\kappa) \To \hm$, which is a section of the
%fibration, and $\Sigma _0 = \varphi(\m (\kappa))$ is transverse to $\xi$. In
%general, $\r ^3$ endowed with the metric
%$$ g_{\kappa, \tau} = \lambda ^2 (dx^2+ dy^2) + (\tau \lambda (y dx -xdy) +dz)^2
%,$$where
%$$ \lambda = \frac{1}{1+\frac{\kappa }{4}(x^2 + y^2)} ,$$is the universal cover of
%$\hm$. Thus, in our case, we make
%
%$$\varphi(\m (\kappa)) = \set{(x,y,z) \in (\r ^3 , g_{\kappa , \tau}) \, : \, z =0
%}.$$
%
%Let $\Omega \subset \m (\kappa)$ be a domain. The Killing graph of a function $u \in
%C^2 (\Omega)$ is the surface
%$$ \Sigma = \set{\gamma _{u(p)}(p) \, ; \, p \in \Sigma _0 }, $$where $u (p) = u (\pi
%(x))$ when $p = \pi (x)$. Henceforth, we identify domains in $\m (\kappa)$ with its
%lift to $\Sigma _0$ and when we say {\it graph} we mean {\it Killing graph}.
%
%So, we have
%
%{\bf Lemma 2.1 \cite{LH}.} {\it Let $\Sigma $ be a surface of $\hm$, transverse to
%the fibers of $\xi$. Let $N$ be a unit normal vector field along $\Sigma$. Then
%\begin{equation}\label{divH}
%{\rm div}_{\kappa}(\pi _\ast N) = 2 H,
%\end{equation}here ${\rm div}_{\kappa}$ stands for the divergence operator in $\m (\kappa)$.}

\begin{definition}
We say that $\Sigma \subset \hm$ is a vertical cylinder over $\alpha$ if $\Sigma =
\pi ^{-1} (\alpha)$, where $\alpha $ is a curve on $\m (\kappa)$.
\end{definition}

It is not hard to verify that if $\alpha $ is a complete curve of geodesic curvature
$2H$ on $\m (\kappa)$, then $\Sigma = \pi ^{-1}(\alpha)$ is complete and has
constant mean curvature $H$. Moreover, these cylinders are characterized by $\nu
\equiv 0$.

We now state two results about the classification of $H-$surfaces. They will be used
in Sections \ref{SKnonnegative} and \ref{SKnonpositive}, but we prove them in
Section \ref{Snuconstant} and Section \ref{Sqconstant} for the sake of clarity. The
first one concerns $H-$surfaces for which the angle function is constant. However, we need to introduce a family of surfaces that appear in the classification:

\begin{definition}\label{DefSurface}
Let $\mathcal S _{\kappa, \tau}$ be a family of complete $H-$surfaces,  in $\hm$, $\kappa <0$, satisfying for any $\Sigma \in \mathcal S _{\kappa , \tau}$:
\begin{itemize}
\item $4H^2 + \kappa < 0$.
\item $q $ vanishes identically on $\Sigma \in \mathcal S _{\kappa , \tau}$, i.e., $\Sigma$ is invariant by a one parameter family of isometries.
\item $0< \nu ^2 < 1$ is constant along $\Sigma$.
\item $K_e = - \tau ^2$ and $K = (\kappa - 4 \tau ^2) \nu ^2 <0$ are constants along $\Sigma$.
\end{itemize}
\end{definition}

An anonymous referee indicated to us the preprint "Hypersurfaces with a parallel higher fundamental form",  by S. Verpoort who observed that we mistakenly omitted the surfaces $\mathcal S _{\kappa, \tau}$  in a first draft of this paper.

\vspace{.1cm}

\begin{theorem}\label{nuconstant}
Let $\Sigma \subset \hm$ be a complete $H-$surface with constant angle function.
Then $\Sigma$ is either a vertical cylinder over a complete curve of curvature $2H$
on $\m (\kappa)$, a slice in $\hr$ or $\sr$, or $\Sigma \in \mathcal S _{\kappa , \tau}$ with $\kappa < 0$.
\end{theorem}

%{\bf Theorem \ref{nuconstant}.} {\it Let $\Sigma \subset \hm$ be a complete
%$H-$surface with constant angle function. Then $\Sigma$ is either a vertical
%cylinder over a complete curve of curvature $2H$ on $\m (\kappa)$, or a slice in
%$\hr$ or $\sr$.}

\vspace{.1cm}

\begin{remark}
Theorem \ref{nuconstant} improves \cite[Lemma 2.3]{ER} for surfaces in $\hr$.
\end{remark}

Of special interest for us are those $H-$surfaces for which the Abresch-Rosenberg
differential is constant.

\vspace{.1cm}

\begin{theorem}\label{qconstant}
Let $\Sigma \subset \hm$ be a complete $H-$surface with $q$ constant.
\begin{itemize}
\item If $q = 0$ on $\Sigma$, then $\Sigma $ is either a slice in $\hr$ or $\sr$ if $H =0 =
\tau$, or $\Sigma $ is invariant by a one-parameter group of isometries of $\hm$.

Moreover, the Gauss curvature of these examples is
\begin{itemize}
\item If $4H^2 + \kappa > 0$, then $K > 0$ they are the rotationally invariant
spheres.
\item If $4H^2 +\kappa =0 $ and $\nu \equiv  0$, then $K \equiv 0 $ and $\Sigma $ is either
a vertical plane in ${\rm Nil}_3$, or a vertical cylinder over a horocycle in $\hr$
or $\widetilde{{\rm PSL}(2,\c)}$.
\item There exists a point with negative Gauss curvature in the remaining cases.
\end{itemize}
%\begin{itemize}
%\item If $H=0=\tau$, $\Sigma$ is a slice in $\hr$ or $\sr$.
%\item If $4H^2 + \kappa > 0$, $\Sigma$ is a rotational embedded sphere $S_H$
%which also implies that $K>0$.
%\item If $4H^2 + \kappa = 0 $ and $\nu \equiv 0$ on $\Sigma$, $\Sigma$ is a vertical
%cylinder over a complete curve of curvature $|\kappa|$. That is, $\Sigma$ is either
%a vertical cylinder over a straight line in ${\rm Nil}_3$, or a vertical cylinder
%over a horocycle in $\hr$ or $\widetilde{{\rm PLS}(2,\c)}$. Moreover, all these
%examples are flat.
%\item If $4H^2 + \kappa \leq  0$ and $\nu$ is not constant, then $\Sigma $ has a point with negative Gauss
%curvature.
%\end{itemize}
\item If $q \neq 0$ on $\Sigma$, then $\Sigma$ is a vertical cylinder over a complete curve of curvature
$2H$ on $\m (\kappa)$.
\end{itemize}
\end{theorem}

%{\bf Theorem \ref{qconstant}.}  {\it Let $\Sigma \subset \hm$ be a complete
%$H-$surface with $q$ constant.
%\begin{itemize}
%\item If $q = 0$ on $\Sigma$, they can be classified as:
%\begin{itemize}
%\item If $H=0=\tau$, $\Sigma$ is a slice in $\hr$ or $\sr$.
%\item If $4H^2 + \kappa > 0$, $\Sigma$ is a rotational embedded sphere $S_H$
%which also implies that $K>0$.
%\item If $4H^2 + \kappa = 0 $ and $\nu \equiv 0$ on $\Sigma$, $\Sigma$ is a vertical
%cylinder over a complete curve of curvature $|\kappa|$. That is, $\Sigma$ is either
%vertical cylinder over a straight line in ${\rm Nil}_3$, or a vertical cylinder over
%a horocycle in $\hr$ or $\widetilde{{\rm PLS}(2,\c)}$. Moreover, all these examples
%are flat.
%\item If $4H^2 + \kappa \leq  0$ and $\nu$ is not constant, then $\Sigma $ has a point with negative Gauss
%curvature.
%\end{itemize}
%\item $If q \neq 0$ on $\Sigma$, then $\Sigma$ is a vertical cylinder over a complete curve of curvature
%$2H$ on $\m (\kappa)$.
%\end{itemize}}

\section{Complete $H-$surfaces $ \Sigma$ with $K \geq 0$}\label{SKnonnegative}

Here we prove

\begin{theorem}\label{Kgeq0}
Let $\Sigma \subset \hm$ be a complete $H-$surface with $K\geq 0$. Then, $\Sigma$ is
either a rotational sphere (in particular, $4H^2+\kappa >0$), or a complete vertical
cylinder over a complete curve of geodesic curvature $2H$ on $\m (\kappa)$.
\end{theorem}
\begin{proof}
The proof goes as follows: First, we prove that $\Sigma$ is a topological sphere or
a complete non-compact parabolic surface. We show that when the surface is a
topological sphere then it is a rotational sphere. If $\Sigma$ is a complete
non-compact parabolic surface, we prove that it is a vertical cylinder by means of
Theorem \ref{qconstant}.

Since $K\geq 0$ and $\Sigma $ is complete, \cite[Lemma 5]{KO} implies that $\Sigma $
is either a sphere or non-compact and parabolic.

If $\Sigma$ is a sphere, then it is a rotational example (see \cite{AR2} or
\cite{AR}). Thus, we can assume that $\Sigma$ is non-compact and parabolic.

We can assume that $q$ does not vanish identically in $\Sigma$. If $q$ does vanish,
then $\Sigma$ is either a vertical cylinder over a straight line in ${\rm Nil}_3$ or
a vertical cylinder over a horocycle in $\hr$ or $\widetilde{{\rm PLS(2,\c)}}$. Note
that we have used here that $K\geq 0$ and  Theorem \ref{qconstant}.

On the one hand, from the Gauss equation \eqref{gauss}
$$ 0 \leq K = K_e + \tau ^2 + (\kappa -4 \tau ^2)\nu ^2  \leq K_e + \tau ^2 + |\kappa - 4\tau
^2|,$$then
\begin{equation}\label{boundaboveK}
H^2-K_e \leq H^2+\tau^2 + |\kappa - 4\tau ^2| .
\end{equation}

On the other hand, using the very definition of $Q \, dz^2 $, \eqref{boundaboveK}
and the inequality $|\xi _1 + \xi _2|^2 \leq 2(|\xi _1|^2 + |\xi|^2)$ for $\xi_1,
\xi _2 \in \c$, we obtain
\begin{equation*}
\begin{split}
\frac{q}{2} & = \frac{2 |Q|^2}{\lambda ^2 } \leq 4(H^2 + \tau ^2)
\frac{4|p|^2}{\lambda ^2} + (\kappa -4\tau ^2)^2\frac{4|A|^4}{\lambda ^2} \\
 &= 4 (H^2 +\tau ^2)(H^2 -K_e) + \frac{(\kappa -4\tau ^2)^2}{4} (1-\nu ^2)^2\\
&\leq 4 (H^2 +\tau ^2)(H^2 -K_e) + \frac{(\kappa -4\tau ^2)^2}{4}  \\
 &\leq 4(H^2 +\tau ^2)(H^2 +\tau ^2 + |\kappa -4\tau ^2|) + \frac{(\kappa -4\tau ^2)^2}{4}.
\end{split}
\end{equation*}

So, from \eqref{lapGauss}, $\Delta \ln q = 4K\geq 0$ and $\ln q$ is a bounded
subharmonic function on a non-compact parabolic surface $\Sigma$ and since the value
$-\infty$ is allowed at isolated points (see \cite{AS}), $q$ is a positive constant
(recall that we are assuming that $q$ does not vanishes identically). Therefore,
Theorem \ref{qconstant} gives the result.

\end{proof}

\section{Complete $H-$surfaces $ \Sigma$ with $K \leq 0$}\label{SKnonpositive}

\begin{theorem}\label{Kleq0}
Let $\Sigma \subset \hm$ be a complete $H-$surface with $K\leq 0$ and $H^2+\tau^2 -
\left|\kappa -4\tau ^2\right| >0$. Then, $\Sigma$ is a complete vertical cylinder
over a complete curve of geodesic curvature $2H$ on $\m (\kappa)$.
\end{theorem}
\begin{proof}
We divide the proof in two cases, $\kappa - 4\tau ^2 < 0$ and $\kappa - 4 \tau
^2>0$.

\vspace{.2cm}

{\bf Case $\kappa - 4 \tau ^2 < 0$:}

\vspace{.2cm}

On the one hand, since $K\leq 0$, we have
$$H^2-K_e \geq H^2+ \tau ^2 + (\kappa -4\tau ^2)\nu ^2 \geq H^2+ \kappa -3\tau ^2,$$
from the Gauss Equation \eqref{gauss}. Therefore, from \eqref{modnuz} and $\kappa -
4\tau ^2 <0$, we obtain:
\begin{equation*}
\begin{split}
q & \geq 4(H^2 + \tau ^2)(H^2 -K_e) + (\kappa -4\tau ^2)(1-\nu ^2 )\left(H^2 +\tau
^2+ H^2 - K_e + \frac{\kappa -4\tau ^2}{4}(1-\nu ^2)\right) \\
 &= (H^2-K_e)\left( 4H^2 + 4\tau ^2 + (\kappa -4\tau ^2)(1-\nu ^2)\right) \\
 & \qquad + (H^2 +\tau ^2)(\kappa -4\tau ^2)(1-\nu ^2 ) + \frac{(\kappa -4\tau ^2)^2}{4}(1-\nu ^2)^2\\
 & \geq (H^2 + \tau^2 + (\kappa -4\tau ^2)\nu ^2)\left( 4H^2 + 4\tau ^2 + (\kappa -4\tau ^2)(1-\nu ^2)\right)\\
 & \qquad + (H^2 +\tau ^2)(\kappa -4\tau ^2)(1-\nu ^2 ) + \frac{(\kappa -4\tau ^2)^2}{4}(1-\nu
 ^2)^2,
\end{split}
\end{equation*}note that the last inequality holds since $4H^2 + 4\tau ^2 + (\kappa -4\tau ^2)(1-\nu ^2)
\geq 4H^2+\kappa > 0$. $4H^2 + \kappa > 0$ follows from
$$0 < 4(H^2 +\tau ^2)- |\kappa - 4\tau ^2| = 4H^2 +\kappa .$$

Set $a:= H^2 +\tau ^2$ and $b:=\kappa -4\tau ^2$. Define the real smooth function
$f:[-1,1] \To \r$ as
\begin{equation}\label{f}
f(x)=(a+bx^2)(4a + b (1-x^2)) + a b (1-x^2) +\frac{b^2}{4}(1-x^2)^2 .
\end{equation}

Note that $q \geq f(\nu)$ on $\Sigma$, $f(\nu)$ is just the last part in the above
inequality involving $q$. It is easy to verify that the only critical point  of $f$
in $(-1,1)$ is $x=0$. Moreover,
$$f(0)=(4a+b)^2/4 > 0 \quad   \mbox{ and } \quad
f(\pm 1) = 4a (a+b) > 0.$$

Actually, $f : \r \To \r $ has two others critical points, $x = \pm
\sqrt{\frac{4a+b}{3|b|}}$, but here, we have used that
$$ \frac{4a + b}{3 |b|} > 1 , $$since $0 < 4(H^2 +\kappa -3\tau ^2) = (4H^2 +\kappa)- 3 |\kappa -4\tau ^2|
= (4a +b) - 3 |b|$.

So, set $c= {\rm min}\set{f(0),f(\pm 1)} >0$, then
$$ q \geq f(\nu) \geq c > 0 .$$

Now, from \eqref{lapGauss} and $q \geq c >0 $ on $\Sigma$, it follows that $ds^2 =
\sqrt{q}I$ is a complete flat metric on $\Sigma$ and
\begin{equation*}
\Delta ^{ds^2} \ln q = \frac{1}{\sqrt{q}} \Delta \ln q = \frac{4K}{\sqrt{q}} \leq 0.
\end{equation*}
%$(\Sigma , ds^2)$ is parabolic. Since parabolicity is preserved by conformal changes
%of the metric, $(\Sigma , I)$ is parabolic. Moreover, from \eqref{lapGauss} and that
%$K\leq 0$, we get $$ \Delta \ln q \leq  K \leq 0 .$$

Since $q$ is bounded below by a positive constant and $(\Sigma , ds^2)$ is
parabolic, then $\ln q $ is  constant which implies that $q$ is a positive constant
(recall $q$ is bounded below by a positive constant). Thus, the result follows from
Theorem \ref{qconstant}. The case $\kappa - 4\tau ^2 < 0$ is proved.

\vspace{.2cm}

{\bf Case $\kappa - 4\tau ^2 >0$: }

\vspace{.2cm}

Set $w_1:= 2(H+i\tau)\dfrac{p}{\lambda}$ and $w_2 : = (\kappa -4\tau
^2)\dfrac{A^2}{\lambda}$, i.e., $ q = 4 |w_1 - w_2|^2$. Then
\begin{equation*}
\begin{split}
|w_1|^2 &= (H^2 +\tau ^2)(H^2 - K_e) \geq (H^2 + \tau ^2)^2 \\[3mm]
|w_2|^2 &= \frac{(\kappa -4\tau ^2)^2}{16}(1-\nu ^2)^2 \leq \left( \dfrac{\kappa
-4\tau^2}{4}\right)^2 ,
\end{split}
\end{equation*}where we have used that $H^2 - K_e \geq H^2 +\tau ^2 + (\kappa -4\tau ^2)\nu
^2\geq H^2 + \tau ^2$, since $K\leq 0$ and $\kappa - 4\tau ^2 >0$.

Let us recall a well known inequality for complex numbers, let $\xi _1 , \xi _2 \in
\c$ then $|\xi _1 + \xi _2|^2 \geq \left||\xi _1|-|\xi _2|\right|^2$. Thus,
\begin{equation*}
\begin{split}
\frac{1}{4} q & \geq \left| |w_1| - |w_2|\right|^2 \geq  \left| (H^2 + \tau ^2) -
\frac{|\kappa -4\tau ^2|}{4} \right|^2 \\
   &= \frac{1}{16}\left| 4(H^2+\tau ^2) - |\kappa -4\tau ^2| \right|^2 >0 .
\end{split}
\end{equation*}

%Now, from \eqref{lapGauss} and $q \geq \frac{1}{4}\left| 4(H^2+\tau ^2) - |\kappa
%-4\tau ^2| \right|^2 >0$ on $\Sigma$, it follows that $ds^2 = \sqrt{q}I$ is a
%complete flat metric on $\Sigma$, then $(\Sigma , ds^2)$ is parabolic. Since
%parabolicity is preserved by conformal changes on the metric, $(\Sigma , I)$ is
%parabolic. Moreover, again from \eqref{lapGauss} and that $K\leq 0$, we get
%$$ \Delta q \leq q K \leq 0 ,$$ since $q$ is bounded below by a positive constant and
%$(\Sigma , I)$ is parabolic, then $q $ is a positive constant. Thus, the result
%follows from Theorem \ref{qconstant}.

So, $q$ is bounded below by a positive constant then, arguing as in the previous
case, $q$ is constant. Thus, the result follows from Theorem \ref{qconstant}. The
case $\kappa - 4\tau ^2 > 0$ is proved.

\end{proof}

\begin{remark}
Note that in the above Theorem, in the case $\kappa - 4\tau^2 > 0$, we only need to
assume $4(H^2+\tau ^2) -|\kappa -4\tau ^2|>0$.
\end{remark}
\section{Complete $H-$surfaces with constant angle function}\label{Snuconstant}

We classify here the complete $H-$surfaces in $\hm$ with constant angle function.
The purpose is to take advantage of this classification result in the next Section.

\vspace{.5cm}

{\bf Theorem \ref{nuconstant}.} {\it Let $\Sigma \subset \hm$ be a complete
$H-$surface with constant angle function. Then $\Sigma$ is either a vertical
cylinder over a complete curve of curvature $2H$ on $\m (\kappa)$, a slice in
$\hr$ or $\sr$, or $\Sigma \in \mathcal S _{\kappa , \tau}$ with $\kappa < 0$ (see Definition \ref{DefSurface}).}

\vspace{.5cm}

\begin{proof}
We can assume that $\nu \leq 0$. We will divide the proof in three cases:

\begin{itemize}
\item $\nu = 0$: In this case, $\Sigma $ must be a vertical cylinder over a complete curve of
geodesic curvature $2H$ on $\m (\kappa)$.

\item $\nu = -1$: From \eqref{c2}, $\tau =0 $ and $H=0$, then $\Sigma $ is a slice in
$\hr$ or $\sr$.

\item $-1 < \nu < 0$: We prove here that $\Sigma \in \mathcal S _{\kappa, \tau}$ with $\kappa <0$. From \eqref{c3}, we
have
\begin{equation}\label{Ap}
(H-i\tau) A = - \frac{2 p}{\lambda} \Ab
\end{equation}then
$$ H^2 +\tau ^2= \frac{4|p|^2}{\lambda ^2} = H^2 -K_e $$since $|A|^2 \neq 0$ from
\eqref{c4}, so $K_e=-\tau ^2$ on $\Sigma$.

Thus, from \eqref{deltanu}, we have
\begin{equation}\label{Hkappa}
4 H^2 + 4\tau ^2+ (\kappa -4\tau^2 ) (1- \nu ^2 ) = 0 .
\end{equation}

Now, using the definition of $q$, \eqref{Ap}, \eqref{Hkappa} and $K_e= -\tau ^2$, we
have
\begin{equation*}
\begin{split}
 q &= \frac{4|Q|^2}{\lambda ^2} = 4(H^2+\tau ^2)\frac{4|p|^2}{\lambda ^2}+
 (\kappa - 4 \tau ^2)^2\frac{4|A|^4}{\lambda ^2}\\
 & \qquad - 4\frac{\kappa - 4 \tau ^2}{\lambda ^2}\left(2\,(H+i \tau)p \Ab ^2 +2\,(H-i \tau)\pb A ^2\right)\\
 &= 4 (H^2 + \tau ^2)(H^2-K_e)+ (\kappa -4 \tau ^2)^2\frac{(1-\nu ^2)^2}{4} + 2 (\kappa - 4\tau
 ^2)(1-\nu ^2)(H^2 + \tau ^2) \\
 &= \frac{1}{4}\left(4H^2 + (\kappa -4\tau ^2)(1-\nu ^2) + 4\tau ^2\right)^2 =
 0
\end{split}
\end{equation*}that is, $q $ vanishes identically on $\Sigma$. Moreover, from \eqref{Hkappa}, we can see that $4H^2 + \kappa <0$, that is, $\kappa <0$. Therefore, $\Sigma \in \mathcal S _{\kappa , \tau}$, $\kappa <0$.

%If $\tau \neq 0$, from \eqref{lapGauss}, $\Sigma $ is flat. But, from the Gauss
%Equation \eqref{gauss} and $K_e=-\tau^2$, $(\kappa - 4\tau ^2)\nu ^2 =0$, which is a
%contradiction. So, $\tau =0 $ and $q\equiv 0$ on $\Sigma$. Recall that we are
%assuming that $0<\nu^2<1$, and moreover $\kappa -4\tau ^2 \neq 0$.
%
%It is well known that a complete $H-$surface with vanishing Abresch-Rosenberg
%differential should be a complete rotational example (see \cite{AR2}). But $\Sigma$
%is also flat, so $\Sigma$ is a vertical cylinder in $\m (\kappa) \times \r$ or the
%complete flat rotational example in $\hr$ described in \cite[Section 3]{AEG2}.
%However, such a complete example has not constant mean curvature. So, in both cases
%we have a contradiction.
\end{itemize}
\end{proof}

\section{Complete $H-$surfaces with $q$ constant}\label{Sqconstant}

Here, we prove the classification result for complete $H-$surfaces in $\hm$ employed
in the proof of Theorem \ref{Kgeq0} and Theorem \ref{Kleq0}.

\vspace{.5cm}

{\bf Theorem \ref{qconstant}.}  {\it Let $\Sigma \subset \hm$ be a complete
$H-$surface with $q$ constant.
\begin{itemize}
\item If $q = 0$ on $\Sigma$, then $\Sigma $ is either a slice in $\hr$ or $\sr$ if $H =0 =
\tau$, or $\Sigma $ is invariant by a one-parameter group of isometries of $\hm$.

Moreover, the Gauss curvature of these examples is
\begin{itemize}
\item If $4H^2 + \kappa > 0$, then $K > 0$ they are the rotationally invariant
spheres.
\item If $4H^2 +\kappa =0 $ and $\nu \equiv  0$, then $K \equiv 0 $ and $\Sigma $ is either
a vertical plane in ${\rm Nil}_3$, or a vertical cylinder over a horocycle in $\hr$
or $\widetilde{{\rm PSL}(2,\c)}$.
\item There exists a point with negative Gauss curvature in the remaining cases.
\end{itemize}
%\begin{itemize}
%\item If $H=0=\tau$, $\Sigma$ is a slice in $\hr$ or $\sr$.
%\item If $4H^2 + \kappa > 0$, $\Sigma$ is a rotational embedded sphere $S_H$
%which also implies that $K>0$.
%\item If $4H^2 + \kappa = 0 $ and $\nu \equiv 0$ on $\Sigma$, $\Sigma$ is a vertical
%cylinder over a complete curve of curvature $|\kappa|$. That is, $\Sigma$ is either
%a vertical cylinder over a straight line in ${\rm Nil}_3$, or a vertical cylinder
%over a horocycle in $\hr$ or $\widetilde{{\rm PLS}(2,\c)}$. Moreover, all these
%examples are flat.
%\item If $4H^2 + \kappa \leq  0$ and $\nu$ is not constant, then $\Sigma $ has a point with negative Gauss
%curvature.
%\end{itemize}
\item If $q \neq 0$ on $\Sigma$, then $\Sigma$ is a vertical cylinder over a complete curve of curvature
$2H$ on $\m (\kappa)$.
\end{itemize}}

\vspace{.5cm}

The case $q=0$ has been treated extensively when the target manifold is a product
space, but is has not been established explicitly  when $\tau \neq 0$. So, we
assemble the results in \cite{AR}, \cite{AR2} for the readers convenience.

\begin{lemma}\label{lemaq=0}
Let $\Sigma \subset \hm $ be a complete $H-$surface whose Abresch-Rosenberg
differential vanishes. Then $\Sigma $ is either a slice in $\hr$ or $\sr$ if $H =0 =
\tau$, or $\Sigma $ is invariant by a one-parameter group of isometries of $\hm$.

Moreover, the Gauss curvature of these examples is
\begin{itemize}
\item If $4H^2 + \kappa > 0$, then $K > 0$ they are the rotationally invariant
spheres.
\item If $4H^2 +\kappa =0 $ and $\nu \equiv  0$, then $K \equiv 0 $ and $\Sigma $ is either
a vertical plane in ${\rm Nil}_3$, or a vertical cylinder over a horocycle in $\hr$
or $\widetilde{{\rm PSL}(2,\c)}$.
\item There exists a point with negative Gauss curvature in the remaining cases.
\end{itemize}
\end{lemma}
\begin{proof}
The idea of the proof for product spaces that we use below, can be found in
\cite{dCF} and \cite{FM}.

If $H=0=\tau$, from the definition of the Abresch-Rosenberg differential, we have
$$ 0 = - (\kappa - 4\tau )A^2 ,$$that is, $\nu ^2 = \pm 1$ using \eqref{c4}. Thus,
$\Sigma $ is a slice in $\hr$ or $\sr$.

If $H \neq 0 $ or $\tau \neq 0$, we have
\begin{equation}\label{aux0}
 2(H +i \tau)p = (\kappa - 4 \tau ^2)A^2 ,
\end{equation}from where we obtain, taking modulus,
\begin{equation}\label{aux1}
H^2 - K_e = \frac{(\kappa -4\tau ^2)^2(1-\nu^2)^2}{16(H^2 +\tau ^2)}
\end{equation}

Replacing \eqref{aux0} in \eqref{c3},
$$ (H+ i \tau )\nu _z = -\frac{1}{4}(4H^2 +\kappa - (\kappa - 4\tau ^2)\nu^2)A
,$$and taking modulus,
\begin{equation}\label{aux2}
|\nu _z|^2 = g(\nu)^2|A|^2, \, \, \,  g(\nu ) = \frac{4H^2 +\kappa -(\kappa - 4\tau
^2)\nu ^2}{4\sqrt{H^2 +\tau ^2}}.
\end{equation}

Assume that $\nu$ is not constant. Let $p\in \Sigma $ be a point where $ \nu _z (p)
\neq 0$ and let $\mathcal{U}$ be a neighborhood of that point $p$ where $\nu _z \neq
0$ (we can assume $\nu ^2 \neq 1 $ at $p$). In particular, $g(\nu ) \neq 0 $ in
$\mathcal{U}$ from \eqref{aux2}. Now, replacing \eqref{aux2} in \eqref{c4}, we
obtain
\begin{equation}\label{aux3}
\lambda = \frac{4 |\nu _z|^2}{(1-\nu ^2)g(\nu)^2} .
\end{equation}

Thus, putting \eqref{aux1} and \eqref{aux3} in the Jacobi equation \eqref{deltanu}
\begin{equation}\label{aux4}
\nu _{z \zb} = -2\frac{\nu |\nu _z|^2}{1-\nu ^2} .
\end{equation}

So, define the real function $s:={\rm arctgh} (\nu)$ on $\mathcal U$. Such a
function is harmonic by means of \eqref{aux4}, thus we can consider a new conformal
parameter $w$ for the first fundamental form so that $s = {\rm Re}(w)$, $w= s + i
t$.

Since $\nu = {\rm tgh}(s)$ by the definition of $s$, we have that $\nu \equiv \nu
(s)$, i.e., it only depends on one parameter. Thus,  we have $\lambda \equiv \lambda
(s)$ and $T \equiv T(s)$ from \eqref{aux3} and \eqref{aux2} respectively, and $p
\equiv p(s)$ by the definition of the Abresch-Rosenberg differential. That is, all
the fundamental data of $\Sigma $ depend only on $s$.

Now, let $\mathcal U $ be a simply connected domain on $\Sigma$ and $\mathcal V
\subset \r ^2$, a simply connected domain of a surface $S$, so that $\psi _0 :
\mathcal V \To \mathcal U \subset \hm $. We parametrice $\mathcal V$ by the
parameters $(s,t)$ obtained above. Then, the fundamental data (see \cite[Theorem
2.3]{FM}) $\set{\lambda _0 , p_0, T_0 , \nu _0}$ of $\psi _0$ are given by
\begin{equation*}
\left\{\begin{matrix}
 \lambda _0 (s,t) &=& \lambda (s) \\
 p_0(s,t) &=& p(s)\\
 T_0(s,t) &=& a(s)\partial _s \\
 \nu _0(s,t) &=& \nu (s) ,
\end{matrix}\right.
\end{equation*}where $a(s) $ is a smooth function.

Set $\bar t \in \r$ and let $\mathbf i _{\bar t}: \r ^2 \To \r ^2$ be the
diffeomorphism given by
$$ \mathbf i _{\bar t}(s,t) : = (s , t + \bar t) , $$and define
$\psi _{\bar t} := \psi _0 \circ \mathbf i _{\bar t}$. Then, the fundamental data
$\set{\lambda _{\bar t} , p_{\bar t}, T_{\bar t} , \nu _{\bar t}}$ of $\psi _{\bar
t} $ are given by
\begin{equation*}
\left\{\begin{matrix}
 \lambda _{\bar t} (s,t) &=& \lambda (s) \\
 p_{\bar t}(s,t) &=& p(s)\\
 T_{\bar t}(s,t) &=& a(s)\partial _s \\
 \nu _{\bar t}(s,t) &=& \nu (s) ,
\end{matrix}\right.
\end{equation*}that is, both fundamental data match at any point $(s,t) \in \mathcal V$.
Therefore, using \cite[Theorem 4.3]{D}, there exists an ambient isometry $\mathcal I
_{\bar t} : \hm \To \hm $ so that
$$ \mathcal I _{\bar t} \circ \psi _0 = \psi_0 \circ \mathbf i _{\bar t}, \text{ for all } \bar t \in \r,$$
thus the surface is invariant by a one parameter group of isometries.

Let us prove the claim about the Gauss curvature. Using the Gauss Equation
\eqref{gauss} in \eqref{aux1}, one gets
$$ H^2 + \tau ^2 + (\kappa - 4\tau^2 )\nu ^2 - K =
\dfrac{(\kappa -4\tau ^2)^2 (1-\nu ^2)^2}{16(H^2+\tau^2)}. $$

Set $a:= 4(H^2 + \tau ^2 )$ and $b:= \kappa - 4\tau ^2$, then one can check easily
that the above equality can be expressed as
\begin{equation}\label{Ch5:Gaussab}
 4a K = a^2 - b^2 + (2a+b)^2 - (2a +b(1-\nu ^2))^2 .
\end{equation}

So, if $4H^2 + \kappa >0$ then $a > |b|$ and $K >0$, that is, $\Sigma$ is a
topological sphere since it is complete. If $4H^2 + \kappa = 0$, $a= - b$ and the
equation reads as
$$ 4a K = a^2(1-(1+\nu ^2)^2),$$that is, $\Sigma $ has a point with negative Gauss
curvature unless $\nu \equiv 0$.

If $4H^2 +\kappa < 0$, one can check that $a^2- b^2 = (a-b)(a+b) <0$ since $a+b >0 $ and $a-b <0$. So, if ${\rm inf}_\Sigma \set{\nu ^2} =0 $ then, from \eqref{Ch5:Gaussab}, $\Sigma $ has a point with negative curvature. Therefore, to finish this lemma, we shall prove that: 
\begin{quote}
{\bf Claim:} There are no complete surfaces with constant mean curvature $4H^2 + \kappa <0$ in $\hm$, $\kappa <0$, with $q \equiv 0$, $K \geq 0$ and ${\rm inf}_\Sigma \set{\nu ^2} =c  >0$.
\end{quote}
{\it Proof of the Claim:} Assume such a surface $\Sigma$ exists. Since we are assuming that $K \geq 0$ and $\Sigma $ is complete, then $\Sigma $ is parabolic and noncompact. If $\Sigma$ were compact we would have a contradiction with the fact that ${\rm inf}_\Sigma \set{\nu ^2} =c>0 $ and $4H^2 + \kappa < 0$. 

Since $q$ vanishes identically on $\Sigma$, ${\rm arctanh} (\nu)$ is a bounded harmonic function on $\Sigma$ and so, $\nu $ is constant. This implies that $K \equiv 0$ and $c< \nu ^2 < 1$ is constant on $\Sigma$. So, the projection $\pi : \Sigma \to \mk $ is a global diffeomorphism and a quasi-isometry. This is impossible since $\Sigma$ is parabolic and $\mk $, $\kappa < 0$, is hyperbolic. Therefore, the Claim is proved and so, the lemma is proved.
\end{proof}

\begin{proof}[Proof of Theorem \ref{qconstant}]

We focus on the case $q \neq 0 $ because Lemma \ref{lemaq=0} gives the
classification when $q=0$.

%When $q=0$, the above classification is a compilation of results in \cite{AR},
%\cite{AR2} and \cite{FM2}. The claim about the Gauss curvature of these surfaces can
%been checked as follows. By the definition of the Abresch-Rosenberg differential,
%$$ Q = 2(H+ i\tau)p - (\kappa -4\tau ^2)A^2 ,$$we first can check that if $H=0$ and
%$\tau =0 $, then $\Sigma$ is a slice. Otherwise, we obtain
%$$ 2(H+i\tau)p = (\kappa - 4 \tau ^2)A^2 ,$$and taking modulus
%$$ (H^2 +\tau ^2)(H^2 - K_e) = \frac{(\kappa -4\tau ^2)^2}{16}(1-\nu ^2)^2 .$$
%
%Using now the Gauss Equation \eqref{gauss},
%$$ (H^2+\tau ^2)(H^2 + \tau ^2 + (\kappa -4\tau ^2)\nu ^2 - K) = \frac{(\kappa -4\tau ^2)^2}{16}(1-\nu ^2)^2 .$$
%
%Set $a:= 4(H^2+\tau ^2)$ and $b:= \kappa - 4\tau ^2$, then we can rewrite the above
%equation as
%$$ 4 a K =  a^2 + 4 ab \nu ^2 - b^2 (1-\nu^2)^2= a^2 - b^2 + (2a +b)^2 - (2a + b(1-\nu^2))^2
%,$$from where we obtain the claim on the Gauss curvature.

Suppose $\nu$ is not constant in $\Sigma$. Since $q = c^2 >0$, we can consider a
conformal parameter $z$ so that $\meta{ \cdot }{ \cdot }=|dz|^2$ and $Q \, dz^2= c
\, dz^2 $ on $\Sigma$. Thus,
$$Q = c = 2(H+i\tau)p - (\kappa - 4\tau ^2)A^2 .$$

First, note that we can assume that $H\neq 0$ or $\tau \neq 0$, otherwise $\nu$
would be constant. So, from \eqref{c3}, we have
$$ (H+i\tau)\nu _z = -(H^2 +\tau ^2 +\frac{\kappa - 4\tau ^2}{4}(1-\nu ^2))A - c \Ab
,$$where we have used $2(H +i \tau)p = c + (\kappa - 4\tau ^2)A^2$. That is,
\begin{equation}\label{nuqconst}
16(H^2 + \tau ^2)\norm{\nabla \nu}^2 = \left( g(\nu)+ 4c \right)^2(1-\nu ^2) ,
\end{equation}where
\begin{equation}\label{gnu}
g(\nu ):= 4H^2+ \kappa - (\kappa -4\tau ^2)\nu ^2 .
\end{equation}

%Combining \eqref{modnuz} and \eqref{nuqconst}, we get that $\nu$ is constant since
%$q$ is constant. Therefore, by means of Theorem \ref{nuconstant}, $\Sigma$ is a
%complete vertical cylinder.

From \eqref{lapGauss}, $\Sigma$ is flat and $H^2-K_e = H^2 + \tau ^2 + (\kappa
-4\tau ^2)\nu ^2$ by \eqref{gauss}, joining this last equation to \eqref{modnuz} we
obtain using the definition of $g(\nu)$ given in \eqref{gnu}
\begin{equation}\label{mod2}
\norm{\nabla \nu}^2 = \dfrac{g(\nu)^2}{4(\kappa -4\tau ^2)}+ \nu ^2 g(\nu)
-\dfrac{c^2}{\kappa -4\tau ^2}.
\end{equation}

Putting together \eqref{nuqconst} and \eqref{mod2} we obtain a polynomial expression
in $\nu ^2$ with coefficients depending on $a:= 4(H^2+ \tau ^2)$, $b:= \kappa -
4\tau ^2$ and $c$,
$$ P(\nu ^2)  := C(a,b,c)\nu ^6 + \text{ lower terms } = 0 ,$$but one can easily
check that the coefficient in $\nu ^6$ is $C(a,b,c)= - a^{-1} b^2 \neq 0$, a
contradiction. Thus $\nu $ is constant, and so, by means of Theorem
\ref{nuconstant}, $\Sigma$ is a vertical cylinder over a complete curve of curvature
$2H$.

\end{proof}

\section{Appendix}

Let $\Sigma $ be a connected Riemannian surface. We establish in this Appendix a
result which we think is of independent interest, concerning differential operators
of the form $\Delta + g$, acting on $C^2(\Sigma)-$functions, where $\Delta $ is the
Laplacian with respect to the Riemannian metric on $\Sigma$ and $g \in C^0
(\Sigma)$.

\begin{lemma}\label{lappendix}
Let $g \in C^{0}(\Sigma)$, $v \in C^{2}(\Sigma )$ such that $\norm{\nabla v}^2 \leq
h \, v^2$ on $\Sigma$, $h$ is a non-negative continuous function on $\Sigma$, and
$\Delta v + g v =0$ in $\Sigma$. Then either $v$ never vanishes or $v$ vanishes
identically on $\Sigma$.
\end{lemma}
\begin{proof}
Set $\Omega = \set{p \in \Sigma \, : \, \, v (p) =0}$. We will show that either
$\Omega = \emptyset$ or $\Omega = \Sigma$.

So, let us assume that $\Omega \neq \emptyset $. If we prove that $\Omega$ is an
open set then, since $\Omega$ is closed and $\Sigma$ is connected, $\Omega =
\Sigma$. Let $p \in \Omega$ and $\mathcal{B}(R) \subset \Sigma$ be the geodesic ball
centered at $p$ of radius $R$. Such a geodesic ball is relatively compact in
$\Sigma$.

Set $\phi = v^2 /2 \geq 0$. Then
$$ \Delta \phi = v \Delta v + \norm{\nabla v}^2 = -g v^2 + \norm{\nabla v}^2\leq - 2(g -h) \phi
,$$that is,
\begin{equation}\label{phi}
-\Delta \phi - 2(g -h) \phi  \geq 0 .
\end{equation}

Define $\beta := {\rm min}\set{{\rm inf}_{\Omega} \set{2(g-h)} , 0} \leq 0$. Then,
$\psi = -\phi $ satisfies
$$ \Delta \psi + \beta \psi = -\Delta \phi - \beta \phi  \geq -\Delta \phi - 2(g-h) \phi \geq 0
,$$where we have used \eqref{phi}.

Since we are assuming that $v $ has a zero at an interior point of $\mathcal{B}(R)$,
$\beta \leq 0 $ and $\psi $ has a non-negative maximum at $p$, the Maximum Principle
\cite[Theorem 3.5]{GT} implies that $\psi $ is constant and so $v$ is constant as
well, i.e, $v \equiv 0$ in $\mathcal{B}(R)$. Then $\mathcal{B}(R)\subset \Omega$,
and $\Omega $ is an open set. Thus $\Omega = \Sigma $.
\end{proof}

\end{document}